\input amstex
\documentstyle{amsppt}
\loadmsbm

\nologo

\TagsOnRight

\NoBlackBoxes

\input epsf
\input supp-pdf

\define\integer{\operatorname{int}}

\define\spn{\operatorname{span}}

\define\dist{\operatorname{dist}}

\define\card{\operatorname{card}}
\def\floor{\mathbin{\hbox{\vrule height1.2ex width0.8pt depth0pt
        \kern-0.8pt \vrule height0.8pt width1.2ex depth0pt}}}

\font\normal=cmss10 scaled 700

\define\sign{\operatorname{sign}}




\define\Normal{\text{\normal N}}

\font\tenscr=callig15 scaled 800
\font\sevenscr=callig15
\font\fivescr=callig15
\skewchar\tenscr='177 \skewchar\sevenscr='177 \skewchar\fivescr='177
\newfam\scrfam \textfont\scrfam=\tenscr \scriptfont\scrfam=\sevenscr
\scriptscriptfont\scrfam=\fivescr

\font\tenscri=suet14
\font\sevenscri=suet14
\font\fivescri=suet14
\skewchar\tenscri='177 \skewchar\sevenscri='177 \skewchar\fivescri='177
\newfam\scrifam \textfont\scrifam=\tenscri \scriptfont\scrifam=\sevenscri
\scriptscriptfont\scrifam=\fivescri

\hsize = 6.5 true in
\vsize = 9 true in


\topmatter
\title
The set of universal interpolating functions is nowhere dense.
\endtitle
\endtopmatter

\centerline{\smc L\. Olsen}
\centerline{Department of Mathematics}
\centerline{University of St\. Andrews}
\centerline{St\. Andrews, Fife KY16 9SS, Scotland}
\centerline{e-mail: {\tt lo\@st-and.ac.uk}}

\bigskip

\centerline{\smc N\. Pugh}
\centerline{Department of Mathematics}
\centerline{University of St\. Andrews}
\centerline{St\. Andrews, Fife KY16 9SS, Scotland}
\centerline{e-mail: {\tt np80\@st-and.ac.uk}}

\bigskip

\centerline{\smc N\. Strout}
\centerline{Department of Mathematics}
\centerline{University of St\. Andrews}
\centerline{St\. Andrews, Fife KY16 9SS, Scotland}
\centerline{e-mail: {\tt ns247\@st-and.ac.uk}}

\bigskip

\footnote""
{
\!\!\!\!\!\!\!\!
}

\topmatter
\abstract
{
In 1998, Benyamini introduced and proved the existence of 
universal interpolating functions.
In the note we prove that the set of
universal interpolating functions is nowhere dense
in the space of continuous functions on $\Bbb R$.
Several extensions and generalisations are also considered.
}
\endabstract
\endtopmatter

\leftheadtext{
L\. Olsen, N\. Pugh \& N\. Strout
}

\rightheadtext{
The set of universal interpolating functions is nowhere dense..
}

\heading{1. Statement of Results.}\endheading

For a subset $X$ of $\ell^{\infty}$, we define the set
$U(X)$ of $X$-universal (interpolating) functions  by
 $$
 \align
 U(X)
&=
 \Big\{
 f\in C(\Bbb R)
 \,\Big|\,
 \text{
 for all sequences $(x_{n})_{n}\in X$ there is  $t\in\Bbb R$
 }\\
 &\qquad\qquad
 \qquad
 \quad\,\text{such that for all integers $n$, we have}\\
  &\qquad\qquad
 \qquad
 \qquad\qquad
 x_{n}=f(t+n)\Big\}\\
  &=
 \Big\{
 f\in C(\Bbb R)
 \,\Big|\,
 \text{
 for all sequences $(x_{n})_{n}\in X$, we have}\\
  &\qquad\qquad
 \qquad
 \qquad\qquad
\exists t\in\Bbb R\,\,:\,\, \sup_{n}|f(t+n)-x_{n}|=0\Big\}
 \,;
 \endalign
$$
here and below $C(\Bbb R)$ denotes the set of real functions on $\Bbb R$.
%
%
%
%
%
%
Note that if $X\subseteq Y\subseteq\ell^{\infty}$, then
 $$
 U(\ell^{\infty})\subseteq U(Y)\subseteq U(X)\,.
 $$
If $X=\ell^{\infty}$, we will simply call $\ell^{\infty}$-universal functions for universal functions.
Universal functions were introduced by Benyamini [Be] in 1998.
It is not clear that universal functions exist.
However, 
somewhat surprisingly
Benyamini [Be]
proved the existence of such functions
as an example of a non-trivial application of the universal surjectivity of the Cantor set.
Universal functions, and generalisations of universal functions, have recently found applications in several parts of mathematics,
including, for example, the study of $\sigma$-ideals of the space of all real valued sequences 
[NaUz],
and  the study of  Abelian group actions [Fr].

The \lq\lq size" of the set $U(\ell^{\infty})$ of universal functions has also recently been studied.
For example,
Bayart \& Quarta [BaQu] have shown that the elements in $U$ are
generic 
from an algebraic viewpoint. More precisely, Bayart \& Quarta [BaQu]
showed that the set $U(\ell^{\infty})$ is maximally lineable, i\.e\.
there is linear subspace $V$ of $C(\Bbb R)$
with
$\dim V=\dim C(\Bbb R)$ such that
 $
 V\setminus\{0\}\subseteq U$.
There are other notions of generic, noticeably, the topological notion of
generic based on Baire category.
In this note we will show that,
while 
the elements in $U(\ell^{\infty})$ are
generic 
from an algebraic viewpoint, they are, nevertheless, 
non-generic
from a topological viewpoint.
More precisely, we show the following slightly more general result, namely, if
$C(\Bbb R)$ is equipped with the metric
$d_{\infty}$ defined by
 $$
 d_{\infty}(f,g)
 =
 \min(\,\|f-g\|_{\infty}\,,\,1\,)
 $$
 for $f,g\in C(\Bbb R)$, and if $X$ is a subset of $(\ell^{\infty})$ satisfying two mild conditions,
 then $U(X)$ is nowhere dense in $(C(\Bbb R),d_{\infty})$, and therefore, in particular, a set of Baire category 1.
 In fact, we will prove a possibly slightly stronger result.
 To state this result, we first define the the set $U_{\Normal}(X)$
of nearly $X$-universal functions by
  $$
 \align
 U_{\Normal}(X)
&=
 \Big\{
 f\in C(\Bbb R)
 \,\Big|\,
 \text{
 for all sequences $(x_{n})_{n}\in X$, we have}\\
  &\qquad\qquad
 \qquad
 \qquad\qquad
 \inf_{t\in\Bbb R}\sup_{n}|f(t+n)-x_{n}|=0\Big\}\,,
 \endalign
$$
 and note that
  $$
  U(X)\subseteq U_{\Normal}(X)\,.
  $$
 We can now state the main result in this note.
 Below we use the following notation, namely, we write
 $\bold 0=(0,0,0.\ldots)$ and $\bold 1=(1,1,1,\ldots)$.
 
 \bigskip
 
 \proclaim{Theorem 1.1}
 Let $X$ is a subset of $\ell^{\infty}$ satisfying the following two conditions:
 \roster
 \item"(1)"
  $
  \bold 0\in X\,;
  $
 \item"(2)"
 $\spn(\bold 1)+X\subseteq X$.
 \endroster
  Then the set  
 $U_{\Normal}(X)$
of nearly $X$-universal functions is 
nowhere dense in $(C(\Bbb R),d_{\infty})$.
In particular, 
 the set  of
 $U(X)$
of $X$-universal functions is 
nowhere dense in $(C(\Bbb R),d_{\infty})$.
\endproclaim

\bigskip

\noindent
It is clear that $X=\ell^{\infty}$ satisfies conditions (1) and (2) in Theorem 1.1.
However, many much smaller subsets of $\ell^{\infty}$ also satisfy conditions (1) and (2).
For example, it is easily seen that the following subsets $X$ of $\ell^{\infty}$ satisfy conditions (1) and (2) in Theorem 1.1:
 $$
 \align
 X
 &=\spn(\bold 1)\,;\\
 X
 &=
 \{\bold x\in\ell^{\infty}
 \,|\,
 \dist(\bold x,\spn(\bold 1))\le \varepsilon\}
 \,\,
 \text{for $0<\varepsilon$}\,;\\
 X
 &=
 \spn(\bold 1)\cup
 \{\bold x\in\ell^{\infty}
 \,|\,
 \varepsilon_{1}\le\dist(\bold x,\spn(\bold 1))\le \varepsilon_{2}\}
 \,\,
 \text{for $0<\varepsilon_{1}\le\varepsilon_{2}$}\,,
 \endalign
 $$
where we have written 
$\dist(\bold x,M)=\inf_{\bold y\in M}\|\bold x-\bold y\|_{\infty}$ for 
$\bold x\in\ell^{\infty}$ and $M\subseteq\ell^{\infty}$.
It is clear that if $X$ satisfies conditions (1) and (2),
then $\spn(\bold 1)\subseteq X$, and $X=\spn(\bold 1)$ is therefore 
the smallest set $X$ that satisfies conditions (1) and (2).
It  follows from this that the largest set $U(X)$ of $X$-universal functions
where $X$ satisfies conditions (1) and (2) is 
$U(\spn(\bold 1))$, 
and
Theorem 1.1 therefore shows that
 the largest set $U(X)$ of $X$-universal functions
where $X$ satisfies conditions (1) and (2) is, in fact, \lq\lq small", namely, nowhere dense.
Finally, since $\ell^{\infty}$ satisfies conditions (1) and (2) in Theorem 1.1, we immediately obtain the following corollary.

\bigskip

\proclaim{Corollary 1.2}
The set  of
 $U_{\Normal}(\ell^{\infty})$
of nearly universal functions is 
nowhere dense in $(C(\Bbb R),d_{\infty})$.
In particular, 
 the set  of
 $U(\ell^{\infty})$
of universal functions is 
nowhere dense in $(C(\Bbb R),d_{\infty})$.
\endproclaim

\bigskip

\noindent
Of course, the first statement in Theorem 1.1
(saying that 
$U_{\Normal}(X)$
 is 
nowhere dense in $(C(\Bbb R),d_{\infty})$)
is only stronger
than the second statement in Theorem 1.1
(saying that 
$U$
 is 
nowhere dense in $(C(\Bbb R),d_{\infty})(X)$)
provided that
$U(X)$ is a proper subset of $U_{\Normal}(X)$.
We have not been able to prove or disprove this,
and we therefore ask the following question.

\bigskip
 
 \proclaim{Question 1.3}
Is $U(X)$ a proper subset of
 $U_{\Normal}(X)$ for $X\subseteq\ell^{\infty}$?
\endproclaim

\bigskip

In addition to the algebraic notion of 
generic called maximally lineable
and the topological
notion of 
generic based on Baire category, there is a third measure theoretic
notion of 
generic called prevalence.
The notion of prevalence
generalizes the 
notion of 
\lq\lq Lebesgue almost all"
to infinite dimensional completely metrizable vector spaces
where
there is no natural analogue of the Lebesgue measure
and was introduced by Hunt, Sauer \& Yorke [HuSaYo] in the 1990's;
we note that
the same ideas were
also
introduced by Christensen [Ch]
in the 1970's in a more general setting 
using different terminology.
Prevalence
is 
defined as follows.
Let $X$ be a vector space
which is completely metrizable.
We now say that a Borel subset $E$ of $X$ is shy in $X$, or that $X\setminus E$ is prevalent in $X$, if the following holds:
there is Borel measure $\Lambda$ on $X$
such that
$\Lambda(X)>0$
and for all $x\in X$ we have
$\Lambda(x+E)=0$.
In
[HuSaYo]
it is shown that
prevalence satisfies all the properties one would expect from a 
generalization
of 
\lq\lq Lebesgue almost all".
We wonder if $U(X)$ is shy in $C(\Bbb R)$ and ask the following question.

\bigskip

\proclaim{Question 1.4}
Is $U(X)$ shy in $C(\Bbb R)$ for $X\subseteq\ell^{\infty}$?
\endproclaim

\bigskip

We close with a question about classes of functions with which $U(\ell^{\infty})$ 
is closed under addition.
For example,
it is clear that if $f$ is a constant function, then 
 $f+U(\ell^{\infty})\subseteq U(\ell^{\infty})$. Hence
 $U(\ell^{\infty})$ 
is closed under addition with respect to the class of constant functions.
 Are there other functions $f\in C(\Bbb R)$ for which
  $f+U(\ell^{\infty})\subseteq U(\ell^{\infty})$?
 In particular, we ask the following question.

  \bigskip

\proclaim{Question 1.5}
Is it true that
$I+U(\ell^{\infty})\subseteq U(\ell^{\infty})$
where $I\in C(\Bbb R)$ denotes the identity function defined by $I(x)=x$ for all $x$?
Is it true that
$\sin+U(\ell^{\infty})\subseteq U(\ell^{\infty})$?
\endproclaim

\bigskip
\bigskip

\heading{2. Proof of Theorem 1.1}\endheading

\bigskip

Let
 $$
 \align
 C_{\text{\rm aff}}(\Bbb R)
 =
 \Big\{
 f\in C(\Bbb R)\,
 &\Big|\,
 \text{$f$ is piecewise affine and has finitely many points}\\
 &\,\,\,\text{of non-differentiability in each compact interval}
 \Big\}
 \endalign
 $$
Also, recall that we write
$\bold 0=(0,0,0.\ldots)$ and $\bold 1=(1,1,1,\ldots)$.

\bigskip

\proclaim{Proposition 2.1}
Let $X$ be a subset of $\ell^{\infty}$. Then the set $U_{\Normal}(X)$ is closed in $(C(\Bbb R),d_{\infty})$.
\endproclaim
\demo{Proof}\newline
\noindent
We must show that the set
  $$
 \align
 C(\Bbb R)\setminus U_{\Normal}(X)
&=
 \Big\{
 f\in C(\Bbb R)
 \,\Big|\,
 \text{
 there is a sequence $(x_{n})_{n}\in X$, such that}\\
  &\qquad\qquad
 \qquad
 \qquad\qquad
 \inf_{t\in\Bbb R}\sup_{n}|f(t+n)-x_{n}|>0\Big\}
 \endalign
$$
is open.
Let $f\in C(\Bbb R)\setminus U_{\Normal}(X)$.
We must now show that there is $r>0$ such that
$B(f,r)\subseteq C(\Bbb R)\setminus U_{\Normal}(X)$.
Since $f\in C(\Bbb R)\setminus U_{\Normal}(X)$,
 there is a sequence $(x_{n})_{n}\in X$, such that
 $\inf_{t\in\Bbb R}\sup_{n}|f(t+n)-x_{n}|>0$.
 For brevity write
 $c=\inf_{t\in\Bbb R}\sup_{n}|f(t+n)-x_{n}|>0$.
Let
 $r=\frac{c}{4}$ and note that $r>0$.
 We now show that
  $
  B(f,r)
  \subseteq
C(\Bbb R)\setminus U_{\Normal}(X)$.
Let $g\in B(f,r)$.
Fix $s\in\Bbb R$ and observe that since 
 $\sup_{n}|f(s+n)-x_{n}|\ge\inf_{t\in\Bbb R}\sup_{n}|f(t+n)-x_{n}|=c>0$, we can choose $n_{s}\in\Bbb N$ such that
 $|f(s+n_{s})-x_{n_{s}}|\ge\frac{c}{2}$. It follows from this that
  $$
  \align
  \tfrac{c}{2}
  &\le
   |f(s+n_{s})-x_{n_{s}}|\\
  &\le
   |f(s+n_{s})-g(s+n_{s})|+|g(s+n_{s})-x_{n_{s}}|\\
   &\le
   d_{\infty}(f,g)+|g(s+n_{s})-x_{n_{s}}|\\
   &<
   r+|g(s+n_{s})-x_{n_{s}}|\\
   &=
   \tfrac{c}{4}+|g(s+n_{s})-x_{n_{s}}|\,,
\endalign
$$
whence
 $
  \tfrac{c}{2}\le|g(s+n_{s})-x_{n_{s}}|
  \le
  \sup_{n}|g(s+n)-x_{n}|$.
  Finally, since $s\in\Bbb R$ was arbitrary, this implies that
  $0< \tfrac{c}{2}
  \le
  \inf_{s\in\Bbb R}
  \sup_{n}|g(s+n)-x_{n}|$, and so $g\in C(\Bbb R)\setminus U_{\Normal}(X)$.
  \hfill$\square$

\bigskip

\proclaim{Lemma 2.2}
We have $\overline{C_{\text{\rm aff}}(\Bbb R)}=C(\Bbb R)$.
\endproclaim
\demo{Proof}\newline
\noindent
This is well-known.
  \hfill$\square$

\bigskip

\proclaim{Lemma 2.3}
Let $f\in C_{\text{\rm aff}}(\Bbb R)$.
\roster
\item"(1)"
The set $\{y\in\Bbb R\,|\,\text{$f^{-1}(y)$ is uncountable}\}$ is countable.
\item"(2)"
The set $\{y\in\Bbb R\,|\,\text{$f^{-1}(y)$ is countable}\}$ is uncountable.
\endroster
\endproclaim
\demo{Proof}\newline
\noindent
(1) Let $M$
be the set of non-differentiability points of $f$.
Since $f\in C_{\text{\rm aff}}(\Bbb R)$,
it follows that there are real numbers $z_{n}$ with $n\in\Bbb Z$ such that
$z_{n}<z_{n+1}$ for all $n$ and
$M=\{z_{n}\,|\,n\in\Bbb Z\}$.
Next, since $f$ is affine on $[z_{n},z_{n+1})$, we conclude that for all $y\in\Bbb R$, we have:
 $$
 \align
 f^{-1}(y)\cap [z_{n},z_{n+1})=\varnothing
\qquad\qquad\qquad
&
\text{if $y\in f([z_{n},z_{n+1}))$;}\\
f^{-1}(y)\cap [z_{n},z_{n+1})\,\,\text{is a singleton}
\qquad
&
\text{if $y\in f([z_{n},z_{n+1}))$ and $f(z_{n})\not= f(z_{n+1})$;}\\
f^{-1}(y)\cap [z_{n},z_{n+1})=[z_{n},z_{n+1})
\qquad\,\,\,
&
\text{if $y\in f([z_{n},z_{n+1}))$ and $f(z_{n})= f(z_{n+1})$, i\.e\.}\\
&
\text{if $y=f(z_{n})$ and $y=f(z_{n+1})$.}\\
\endalign
$$
It follows immediately from the above that
 $
\{y\in\Bbb R\,|\,\text{$f^{-1}(y)$ is uncountable}\}
\subseteq
\{f(z_{n}\,|\,n\in\Bbb Z\}$, and this implies that
$\{y\in\Bbb R\,|\,\text{$f^{-1}(y)$ is uncountable}\}$ is countable.

\noindent
(2) This follows immediately from (1).
  \hfill$\square$

\bigskip

\proclaim{Lemma 2.4}
Let $X$ be a subset of $\ell^{\infty}$ 
satisfying condition (2) in  Theorem 1.1.
If $f\in C(\Bbb R)\setminus U_{\Normal}(X)$ and $y\in\Bbb R$, then
$f+y\in C(\Bbb R)\setminus U_{\Normal}(X)$.
\endproclaim
\demo{Proof}\newline
\noindent
This follows immediately from the definitions.
  \hfill$\square$ 

\bigskip

\proclaim{Proposition 2.5}
Let $X$ be a subset of $\ell^{\infty}$ 
satisfying conditions (1) and (2) in  Theorem 1.1.
Then we have $C_{\text{\rm aff}}(\Bbb R)
\subseteq \overline{C(\Bbb R)\setminus U_{\Normal}(X)}$.
In particular, we have
$\overline{C_{\text{\rm aff}}(\Bbb R)}
\subseteq \overline{C(\Bbb R)\setminus U_{\Normal}(X)}$.
\endproclaim
\demo{Proof}\newline
\noindent
Let $f\in C_{\text{\rm aff}}(\Bbb R)$.
We must now find a sequence $(f_{n})_{n}$ in $C(\Bbb R)\setminus U_{\Normal}(X)$
such that
$f_{n}\to f$ with respect to $d_{\infty}$.

It follows from Lemma 2.3 that there is a real number $y\in\Bbb R$ such that
$f^{-1}(y)$ is countable. Now put $g=f-y$ and note that
since 
$g^{-1}(0)=f^{-1}(y)$, we deduce that $g^{-1}(0)$ is countable.
For an integer $k$ write
 $$
 \align
 P_{k}
&=
 g^{-1}(0)\cap [5k,5k+1)\,,\\
 Q_{k}
&=
 P_{k}+1\,,\\
 R_{k}
&=
 P_{k}\cup Q_{k}\,.
 \endalign
 $$
Next, note that
 since
 $g\in C_{\text{\rm aff}}(\Bbb R)$ and $g^{-1}(0)$ is countable,
the set $g^{-1}(0)\cap [n,n+1)$ is finite for all $n\in\Bbb Z$.
 In particular, this implies that we can choose
 $s_{k}>0$ such that:
  $$
  \align
 &\text{
  The family $(\,(x-s_{k},x+s_{k})\,)_{x\in R_{k}}$ consists of pairwise disjoint sets;
  }\tag2.1\\
 &\text{
 If $z\in P_{k}$, then $(z-s_{k},z+s_{k})\cap P_{k}=\{z\}$;
  }\tag2.2\\
 &\text{
 If $w\in Q_{k}$, then $(w-s_{k},w+s_{k})\cap Q_{k}=\{w\}$
 and
 $(w-s_{k},w+s_{k})\cap (g^{-1}(0)\setminus Q_{k})=\varnothing$.}\tag2.3
  \endalign
  $$
Now let $r_{k}=\min(\frac{1}{4},\frac{s_{k}}{4})$, and define
 $$
 g_{n,k}
 :
 [2k-r_{k},2k+2+r_{k}]\setminus
 \Bigg(
 \bigcup_{x\in R_{k}}(x-r_{k},x+r_{k})
 \Bigg)
 \to
 \Bbb R
 $$
 by
 $$
 g_{n,k}(x)=\sign(g(x))\,\tfrac{1}{n}\,.
 $$
Next, define
$h_{n,k}:\Bbb R\to\Bbb R$ by
 $$
h_{n,k}(x)
=
\cases
g_{n,k}(x)
&\quad
\text{for}\,\,
x\in
[2k-r_{k},2k+2+r_{k}]\setminus
 (
 \bigcup_{x\in R_{k}}(x-r_{k},x+r_{k})
 )\,;\\
 g_{n,k}(z-r_{k})
 &\quad       
\text{if there is $z\in P_{k}$ such that $x\in[z-r_{k},z]$}\,;\\
 g_{n,k}(z+r_{k})
 &\quad       
\text{if there is $w\in Q_{k}$ such that $x\in[w,w+r_{k}]$}\,;\\
\frac{g_{n,k}(z+r_{k})-g_{n,k}(z-r_{k})}{r_{k}}x+g_{n,k}(z+r_{k})
 &\quad       
\text{if there is $z\in P_{k}$ such that $x\in[z,z+r_{k}]$}\,;\\
\frac{g_{n,k}(w+r_{k})-g_{n,k}(w-r_{k})}{r_{k}}x+g_{n,k}(w-r_{k})
 &\quad       
\text{if there is $w\in Q_{k}$ such that $x\in[w-r_{k},w]$}\,;\\
\frac{g_{n,k}(5k-r_{k})}{r_{k}}x+g_{n,k}(5k-r_{k})
 &\quad       
\text{if  $x\in[5k-2r_{k},5k-r_{k}]$}\,;\\
\frac{g_{n,k}(5k+2+r_{k})}{r_{k}}x
 &\quad       
\text{if  $x\in[5k+2r_{k},5k+2+2r_{k}]$}\,;\\
0
&\quad
\text{otherwise.}
\endcases
$$ 
The somewhat opaque-looking function $h_{n,k}$, is a continuous piecewise linear function constructed so that at points $x\in [5k, 5k+1)$ with $g(x)=0$,
the function $h_{n,k}$ takes value $\pm \frac{1}{n}$, and at points $x\in [5k+1, 5k+2)$ with $g(x-1)=0$, 
the function $h_{n,k}$ takes values $\mp \frac{1}{n}$. 
In Figure 1 below we sketch the graph of $h_{n,k}$ for 
$\card P_{k}=\card Q_{k}=2$.

\midinsert
\centerline
{
\hbox
{
\hskip 0mm
\pdfimage width 15cm{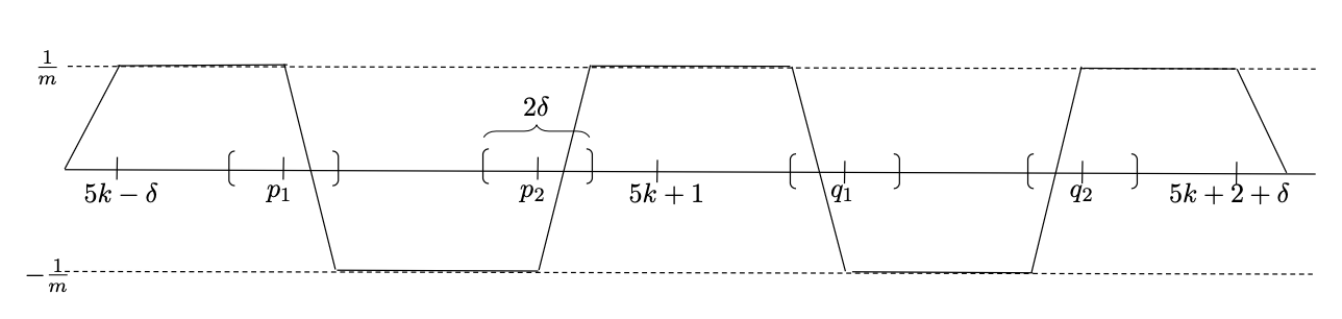}
}
}


\botcaption{\eightpoint \bf Figure 1}
\eightpoint 
A sketch of the graph of $h_{n,k}$ for 
$P_{k}=\{z_{1},z_{2}\}$ and 
$Q_{k}=\{w_{1},w_{2}\}$.
\endcaption
\endinsert

\bigskip

For a fixed $n$, 
the supports of the functions $(h_{n,k})_{k\in\Bbb Z}$ are pairwise disjoint,
and the series $\sum_{k\in\Bbb Z}h_{n,k}$ is therefore well-defined and we can therefore define the function $h_{n}$ by
 $$
 h_{n}=\sum_{k\in\Bbb Z}h_{n,k}\,.
 $$
Finally, let
 $$
 f_{n}=f+h_{n}\,.
 $$
We will now prove that $f_{n}\in C(\Bbb R)\setminus U_{\Normal}(X)$ and that $f_{n}\to f$ with respect to $d_{\infty}$;
this is the contents of the  claims below.

Let
 $$
 g_{n}=g+h_{n}\,.
 $$
 
\bigskip
\noindent
{\it Claim 1. $g_{n}\in C(\Bbb R)\setminus U_{\Normal}(X)$} 

\noindent
{\it Proof of Claim 1.}
In order to prove that $g_{n}\in C(\Bbb R)\setminus U_{\Normal}(X)$,
we must show that there is a sequence $(x_{m})_{m}\in X$ such that
$$
\inf_{t\in\Bbb R}\sup_{m}|g_{n}(t+m)-x_{m}|>0\,.
\tag2.4
$$
Below we will show that the sequence $\bold 0\in X$  satisfies
(2.4), more precisely, we will show that
$$
\inf_{t\in\Bbb R}\sup_{m}|g_{n}(t+m)|\ge\tfrac{1}{n}\,.
\tag2.5
$$
 Assume, in order to reach a contradiction, that (2.5) is not satisfied.
 We can therefore find $t_{0}\in\Bbb R$ such that
 $$
 \sup_{m}|g_{n}(t_{0}+m)|<\tfrac{1}{n}\,.
 \tag2.6
$$

Next choose 
an integer $k_{0}$ such that $t_{0}<5k_{0}$ and let
$m_{0}=5k_{0}-\integer(t_{0})$
where $\integer(x)$ denotes the integer part of $x\in\Bbb R$.
Then $m_{0}=5k_{0}-\integer(t_{0})\in\Bbb N$
and
$t_{0}+m_{0}
=
t_{0}+5k_{0}-\integer(t_{0})
=
t_{0}-\integer(t_{0})+5k_{0}
\in
[5k_{0},5k_{0}+1)$.
This and (2.6) now imply that
 $$
  |g(t_{0}+m_{0})+h_{n}(t_{0}+m_{0})|
  =
 |g_{n}(t_{0}+m_{0})|<\tfrac{1}{n}\,.
 \tag2.7
 $$
Also observe that it follows from the definition of $h_{n}$ that
 $
 \sign(g(x))
 =
 \sign(h_{n}(x))
 $ 
 for all
 $x\in
 [5k_{0},5k_{0}+1)
 \setminus
 (
 \cup_{z\in\ P_{k_{0}}}(z,z+r_{k})
 )
 $, whence
 $$
  |h_{n}(x)|
  \le
 |g(x)+h_{n}(x)|
 \,\,\text{for all}\,\,
 x\in
 [5k_{0},5k_{0}+1)
 \setminus
 \Bigg(
 \bigcup_{z\in P_{k_{0}}}(z,z+r_{k_{0}})
 \Bigg)\,.
 \tag2.8
 $$
On the other hand, the definition of $h_{n}$ also shows that
 $$
 \tfrac{1}{n}
 =
  |h_{n}(x)|
 \,\,\text{for all}\,\,
 x\in
 [5k_{0},5k_{0}+1)
 \setminus
 \Bigg(
 \bigcup_{z\in P_{k_{0}}}(z,z+r_{k_{0}})
 \Bigg)\,.
 \tag2.9
 $$
Combining (2.7) and (2.8), we now deduce that
$$
  \tfrac{1}{n}
  \le
 |g(x)+h_{n}(x)|
 \,\,\text{for all}\,\,
 x\in
 [5k_{0},5k_{0}+1)
 \setminus
 \Bigg(
 \bigcup_{z\in P_{k_{0}}}(z,z+r_{k_{0}})
 \Bigg)\,.
 \tag2.10
 $$
Recalling that
$t_{0}+m_{0}
\in
[5k_{0},5k_{0}+1)$, it follows from (2.7) and (2.10) that
$t_{0}+m_{0}\in \cup_{z\in P_{k_{0}}}(z,z+r_{k_{0}})$, and can therefore find $z_{0}\in P_{k_{0}}$ with
$t_{0}+m_{0}\in (z_{0},z_{0}+r_{k_{0}})$.

Since
$t_{0}+m_{0}\in (z_{0},z_{0}+r_{k_{0}})$ with
$z_{0}\in P_{k_{0}}$, we conclude that
$t_{0}+m_{0}+1\in (w_{0},w_{0}+r_{k_{0}})$ for some
$w_{0}\in Q_{k_{0}}$.
Hence
 $$ 
 \align
 \tfrac{1}{n}
 &>
 \inf_{m}|g_{n}(t_{0}+m)|
 \qquad\qquad
 \qquad\qquad
 \qquad\qquad\,\,\,
 \text{[by (2.6)]}\\
 &\ge
 |g_{n}(t_{0}+m_{0}+1)|\\
 &=
 |g(t_{0}+m_{0}+1)+h_{n}(t_{0}+m_{0}+1)|\\
  &=
 |g(t_{0}+m_{0}+1)+h_{n,k_{0}}(t_{0}+m_{0}+1)|
\qquad\quad
 \text{[since $h_{n}(t_{0}+m_{0}+1)=h_{n,k_{0}}(t_{0}+m_{0}+1)$}\\
 &\quad
\qquad\qquad
\qquad\qquad
\qquad\qquad
\qquad\qquad
\qquad\quad\,\,\,\,
 \text{because}\\
 &\quad
\qquad\qquad
\qquad\qquad
\qquad\qquad
\qquad\qquad
\qquad\quad\,\,\,\, 
 \text{
 $t_{0}+m_{0}+1\in (w_{0},w_{0}+r_{k_{0}})\subseteq  
 [2k_{0}-r_{k_{0}},2k_{0}+2+r_{k_{0}}]$]}\\ 
   &=
 |g(t_{0}+m_{0}+1)+\sign(g(t_{0}+m_{0}+1))\,\tfrac{1}{n}|
\quad\,\,
 \text{[since $h_{n,k_{0}}(t_{0}+m_{0}+1)=\sign(g(t_{0}+m_{0}+1))\,\tfrac{1}{n}$}\\
 &\quad
\qquad\qquad
\qquad\qquad
\qquad\qquad
\qquad\qquad
\qquad\quad\,\,\,\,
 \text{because}\\
 &\quad
\qquad\qquad
\qquad\qquad
\qquad\qquad
\qquad\qquad
\qquad\quad\,\,\,\, 
 \text{
 $t_{0}+m_{0}+1\in (w_{0},w_{0}+r_{k_{0}})\subseteq 
 [2k_{0}-r_{k_{0}},2k_{0}+2+r_{k_{0}}]$]}\\ 
 &\ge
 \tfrac{1}{n}\,.
 \endalign
 $$
This provides the required contradiction.
This completes the proof of Claim 1.
 
 \bigskip
\noindent
{\it Claim 2. $f_{n}\in C(\Bbb R)\setminus U_{\Normal}(X)$} 

\noindent
{\it Proof of Claim 2.}
Since $f_{n}=g_{n}+y$ and since it follows from Claim 1 that $g_{n}\in C(\Bbb R)\setminus U_{\Normal}(X)$, 
we 
now conclude from Lemma 2.4 that
$f_{n}\in C(\Bbb R)\setminus U_{\Normal}(X)$.
This completes the proof of Claim 2.

 \bigskip
\noindent
{\it Claim 3. $f_{n}\to f$ with respect to $d_{\infty}$} 

\noindent
{\it Proof of Claim 3.}
This follows immediately from the fact that
$d_{\infty}(f_{n},f)
=
\|h_{n}\|_{\infty}
=
\frac{1}{n}$.
This completes the proof of Claim 3.

\bigskip

The proof of Proposition 2.5 now follows from Claim 2 and Claim 3.
\hfill$\square$

\bigskip

\demo{Proof of Theorem 1.1}\newline
\noindent
It follows immediately from Lemma 2.2 and Proposition 2.5 that
 $
 C(\Bbb R)
 =
 \overline{C_{\text{\rm aff}}(\Bbb R)}
 \subseteq
 \overline{C(\Bbb R)\setminus U_{\Normal}(X)}$,
i\.e\. $C(\Bbb R)\setminus U_{\Normal}(X)$ is dense, and since it also follows from Proposition 2.1 that
  $U_{\Normal}(X)$ is closed, we now conclude that
  $U_{\Normal}(X)$ is meagre.
  \hfill$\square$

\end